# A Non-Linear Difference Equation for Calculation of the Zeros of the Riemann Zeta-Function on the Critical Line

G. B. da Silva    and    R. V. Ramos
george_barbosa@fisica.ufc.br    rubens.ramos@ufc.br

*Lab. of Quantum Information Technology, Department of Teleinformatic Engineering – Federal University of Ceara - DETI/UFC, C.P. 6007 – Campus do Pici - 60455-970 Fortaleza-Ce, Brazil.*

*Abstract*

In this work, we present a non-linear difference equation for calculation of the zeros of the Riemann's zeta-function on the critical line. Our proposed non-linear map uses the Lambert *W* function and it can be easily implemented in a mathematical software. In order to check the quality of the zeros calculated, we show the factorization of an integer number by the calculation of the discrete cosine Riemann transform.

## 1. Introduction

The Riemann zeta-function has a strong connection with prime numbers and, because of this, it plays an important role in mathematics, physics and computer science [1-6]. The Riemann zeta-function is given by

$$\zeta(s) = \sum_{n=1}^{\infty} \frac{1}{n^s} = \prod_{p} \left(1 - \frac{1}{p^s}\right)^{-1}. \tag{1}$$

In (1) $s = \sigma + it$ is a complex variable and $p$ runs over all prime numbers. Eq. (1) converges only for $\sigma > 1$. The analytical continuation

$$\zeta(s) = \left(1 - 2^{1-s}\right)^{-1} \sum_{n=1}^{\infty} \frac{(-1)^{n+1}}{n^s} \tag{2}$$

permits to calculate $\zeta(s)$ in the region [0,1). Note there is a single pole at $s = 1$. Among the different applications of the zeta-function, one of the most important is the calculation of the prime counting function given by

$$\pi(x) = \sum_{n=1}^{N} \frac{\mu(n)}{n} \left[ Li\left(x^{1/n}\right) - \sum_{\rho} Li\left(x^{\rho/n}\right) + \int_{x^{1/n}}^{\infty} \frac{1}{u(u^2-1)\ln(u)} du - \ln(2) \right]. \tag{3}$$

In (3) $Li(x)$ is the logarithmic integral and $\mu(n)$ is the Möbius function. Furthermore, $\rho$ runs over all non-trivial zeros of the zeta-function (the trivial zeros are the negative even integer numbers). Due to its importance, different numerical methods for calculation of those zeros were developed [7,8]. In particular, in [9] França and LeClair provided an exact equation that should be obeyed by all zeros on the critical line $\sigma = ½$. Besides that, they also showed an analytical equation based on the Lambert $W$ function that could provide a reasonable estimation of the zeros on the critical line. This is estimation is not good enough to reproduce the statistics of the distance between consecutives zeros (GUE statistics) nor to calculate the prime counting function. However, the estimated value could be used as start point of a numerical algorithm for solving the exact equation proposed.

In this direction, the present work proposes a non-linear difference equation also based on the Lambert $W$ function, which provides a very good estimation of the zeros on the critical line. In order to check the quality of the zeros calculated, we used them in the discrete cosine Riemann transform for factorization of an integer number. At last, it is also shown that our non-linear map has a rich dynamic reaching the chaotic behaviour.

This work is outlined as follows: In Section 2, the proposed non-linear map for calculation of the zeros on the critical line is presented. In Section 3, the discrete Riemann cosine transform is introduced and calculated using the zeros provided by our non-linear map. At last, the conclusions are drawn in Section 4.

## 2. The Difference Non-Linear Equation for Calculation of the Zeros on the Critical Line.

There are different strategies for numerical calculations of the zeros on the critical line [7,8]. According to [9] the $n$-th non-trivial zero with positive imaginary part, $½ + it_n$, must obey the following equation

$$\arg \Gamma\left(\frac{1}{4} + \frac{i}{2}t_n\right) - t_n \log \sqrt{\pi} + \lim_{\delta \to 0^+} \arg \zeta\left(\frac{1}{2} + \delta + it_n\right) = \left(n - \frac{3}{2}\right)\pi. \tag{4}$$

In (4), $\Gamma(z)$ is the gamma function that can be written in the Weierstrass' representation as

$$\Gamma(z) = \frac{e^{-\gamma z}}{z} \prod_{k=1}^{\infty} \left(1 + \frac{z}{k}\right)^{-1} e^{z/k}. \tag{5}$$

Asymptotically, when $z$ goes to infinite, the gamma function is given by Stirling's formula

$$\Gamma(z+1) \sim \sqrt{2\pi z} \left(\frac{z}{e}\right)^z. \tag{6}$$

Using (6) in (4), one has the asymptotic form of (4)

$$\frac{t_n}{2\pi}\log\left(\frac{t_n}{2\pi e}\right) + \lim_{\delta \to 0^+} \frac{1}{\pi}\arg\zeta\left(\frac{1}{2}+\delta+it_n\right) = n - \frac{11}{8} \quad (n=1,2,...), \tag{7}$$

that provides an estimation of the zeros, the larger the value of $n$ the better is the estimation. At last, when the term with limit in (7) is not considered, one gets

$$\frac{\hat{t}_n}{2\pi}\log\left(\frac{\hat{t}_n}{2\pi e}\right) = \left(n - \frac{11}{8}\right) \tag{8}$$

whose solutions $\hat{t}_n$ are reasonable estimations of $t_n$. Equation (8) can be solved by using the Lambert $W$ function:

$$\hat{t}_n = 2\pi\left(n-\frac{11}{8}\right)\bigg/ W_0\left(\frac{1}{e}\left(n-\frac{11}{8}\right)\right). \tag{9}$$

Figure 1 shows the first 160,000 zeros obtained from [10], $t_n$, and using (9), $\hat{t}_n$, while Fig. 2 shows the difference between them.

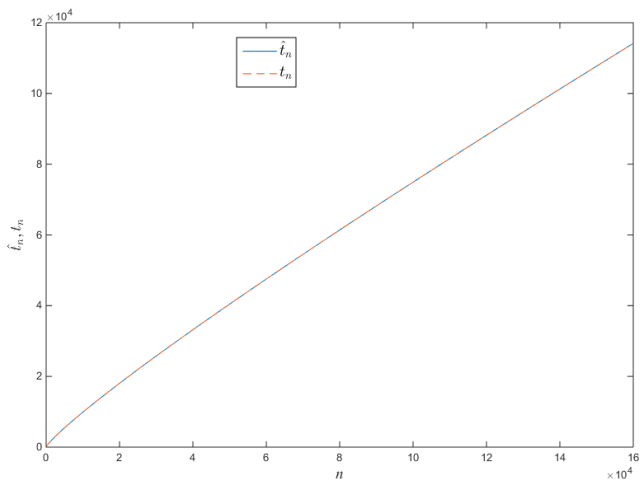

Figure 1 - $t_n$ and $\hat{t}_n$ versus $n$ for the first 160,000 zeros on the critical line.

Observing Figs. (1) and (2) one may think that (9) provides a good estimation, however, this is not true since the prime counting function nor the GUE statistics cannot be correctly obtained by using (9). Nevertheless, the values obtained with (9) can be used as start point in the algorithms that solve (4) and (7).

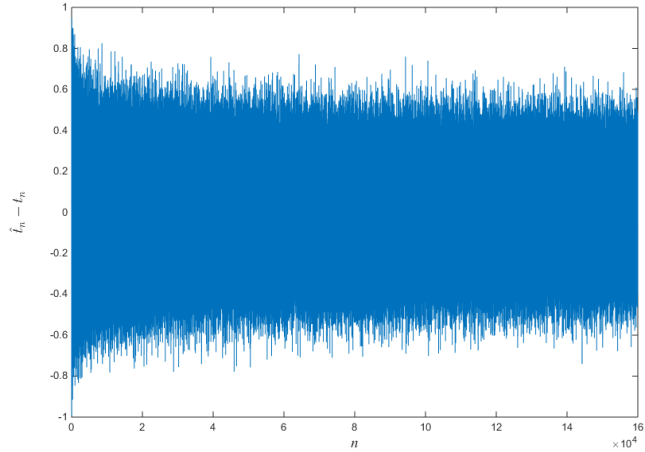

Figure 2 - $\hat{t}_n$- $t_n$ versus $n$ for the first 160,000 zeros on the critical line.

In order to provide a non-linear difference equation that converges to the zeros, we use (7) as inspiration. Consider the following equation

$$\frac{t_n}{2\pi}\log\left(\frac{t_n}{2\pi e}\right)=n-\frac{11}{8}-\frac{\delta}{\pi}\arg\zeta\left(\frac{1}{2}+it_n\right), \qquad (10)$$

whose solution using the Lambert $W$ function is

$$t_n^k = \frac{2\pi\left(n-\frac{11}{8}-\frac{\delta}{\pi}\arg\zeta\left(\frac{1}{2}+it_n^{k-1}\right)\right)}{W_0\left(\frac{1}{e}\left(n-\frac{11}{8}-\frac{\delta}{\pi}\arg\zeta\left(\frac{1}{2}+it_n^{k-1}\right)\right)\right)}. \qquad (11)$$

One can immediately note that we transformed the solution of (10) in a non-linear map. The goal is to obtain better values for $t_n$ by iterating (11). The optimal value of $\delta$ is crucial in order to have $t_n^k$ converging to the correct value $t_n$. One can see in Table 1 the values of $t_n$ [10], $t_n^{k=20}$ (eq. (11)) and $\hat{t}_n$ (eq. (9)) for $n \in \{1, 10{,}000, 50{,}000, 100{,}000\}$, as well the values of $\delta$ used.

Table 1 – Values of $t_n$ [10], $t_n^{k=20}$ and $\hat{t}_n$ for $n \in \{1, 10.000, 50.000, 100.000\}$.

| $n$ | $\delta$ | $t_n$ | $t_n^{k=20}$ | $\hat{t}_n$ |
|---|---|---|---|---|
| 1 | 0.0921796 | 14.134725142000001 | 14.134725496347967 | 14.521346953065633 |
| 10,000 | 0.2639143 | 9877.782654004001 | 9877.782653979717 | 9877.629616492992 |
| 50,000 | 0.1572079 | 40433.68738546200 | 40433.68738541853 | 40433.62056224795 |
| 100,000 | 0.1595388 | 74920.82749899400 | 74920.82749899139 | 74920.89103264698 |

The non-linear map (11) was iterated only 20 times and for all cases considered the initial value used was $t_n^{k=0} = 1$. It can be clearly seen in Table 1 that (11) provides much better estimations of the zeros than (9).

## 3. The Discrete Cosine Riemann Transform

The discrete cosine Riemann transform (DCRT) of a function $x(t)$ is defined as

$$X(k) = \sum_{n=1}^{\infty} x(t_n) \cos\left[\log(k) t_n\right], \tag{12}$$

where $t_n$ are the positive imaginary parts of the non-trivial zeros of the Riemann zeta function. One can note that two different functions can have the same DCRT only if they have the same values at the sample points ($t_n$'s): If $X(k)$ is the DCRT of $x(t)$ and $y(t)$, then

$$\sum_{n=1}^{\infty} x(t_n) \cos\left[\log(k) t_n\right] - \sum_{n=1}^{\infty} y(t_n) \cos\left[\log(k) t_n\right] = \sum_{n=1}^{\infty} \left[x(t_n) - y(t_n)\right] \cos\left[\log(k) t_n\right] = 0. \tag{13}$$

Equation (13) is true for all values of $k$ only if $x(t_n) = y(t_n)$ for all $n$. Now, let us assume that $x(t)$ is given by

$$x(t) = \cos\left[\log(R) t\right]. \tag{14}$$

In this case, the DCRT of $x(t)$ is

$$X(k) = \sum_{n=1}^{\infty} \cos\left[\log(R) t_n\right] \cos\left[\log(k) t_n\right] = \frac{1}{2} \sum_{n=1}^{\infty} \cos\left[\log(Rk) t_n\right] + \frac{1}{2} \sum_{n=1}^{\infty} \cos\left[\log\left(\frac{R}{k}\right) t_n\right]. \tag{15}$$

Observing (15) one can see that $X(k)$ goes to infinite when $R \times k$ or $R/k$ is prime, a prime power or equal to 1, and zero otherwise. For example: 1) if $R = p$, a prime number, $X(k)$ is zero for all values of $k$ except for spikes at $k = p^0, p^1, p^2, p^3, \ldots$ 2) if $R = p_1 p_2$, where $p_1$ and $p_2$ are prime numbers, $X(k)$ is zero for all values of $k$ except for spikes at $k = p_1$, $k = p_2$ and $k = p_1 p_2$. 3) if $R = p_1 p_2 p_3$, where $p_1$, $p_2$ and $p_3$ are prime numbers, $X(k)$ is zero for all values of $k$ except for spikes at $k = p_1 p_2$, $k = p_1 p_3$, $k = p_2 p_3$ and $k = p_1 p_2 p_3$. Hence, the DCRT can be used for integer factorization.

In order to compare the quality of the zeros obtained by (9) and (11) we calculated the DCRT of $x(t) = \cos[\log(131 \times 761) t]$ using only the first 160,000 zeros on the critical line. In Fig. 3 one can see $X(t_n)$ obtained using the zeros from [10], $X(\hat{t}_n)$ calculated with the zeros obtained from (9), and $X(\bar{t}_n)$ where $\bar{t}_n = t_n + \varepsilon$, in which $\varepsilon$ is a uniformly distributed random variable in the interval [-1,1].

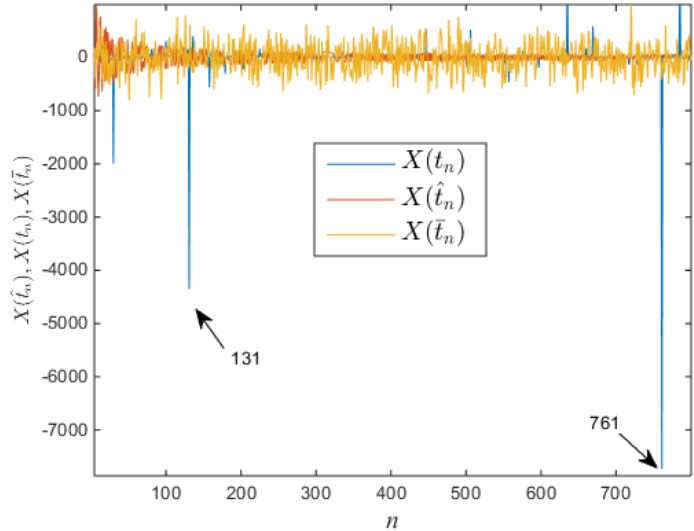

Figure 3 – Discrete cosine Riemann transform of $x(t) = \cos[\log(131 \times 761)t]$. $X(t_n)$ uses the zeros from [10], $X(\hat{t}_n)$ uses the zeros obtained with (9), and $X(\bar{t}_n)$ uses $\bar{t}_n = t_n + \varepsilon$, where $\varepsilon$ is a random variable uniformly distributed in the interval [-1,1].

As one can see in Fig. 3, $\bar{t}_n$ and $\hat{t}_n$ are not good estimations for the correct calculation of the DCRT. They do not show the peaks at $k = 131$ and $k = 761$. Figure 4, shows the DCRT of $x(t) = \cos[\log(131 \times 761)t]$ using only the first 160,000 zeros on the critical line obtained using (11), $X(t_n^{k=20})$. For each zero, (11) is iterated only 20 times with the best value of $\delta$ used. Only the first five digits after the decimal point were considered.

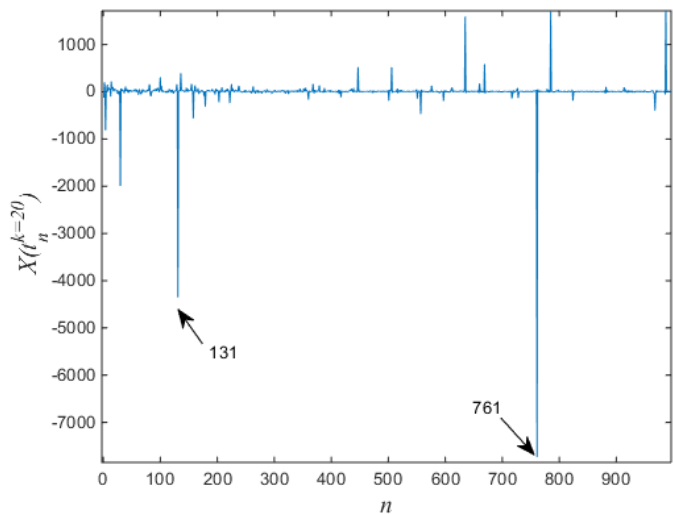

Figure 4 – Discrete cosine Riemann transform of $x(t) = \cos[\log(131 \times 761)t]$, using only the first 160,000 zeros on the critical line using the zeros obtained with (11).

As one can see in Fig. 4, if the optimal $\delta$ value is used, few iterations are required in order to have an estimation of zeros good enough for calculation of the DCRT.

## 4. Non-linear Dynamic

One may also see that (11) is a non-linear map that owns a rich dynamic. Depending on the value of the parameter $\delta$ it can show a stable, periodic or chaotic behaviour. Figures 5 ($\delta = 0.1595388$), 6 ($\delta = 0.5$) and 7 ($\delta = 2$) show for $n = 100{,}000$, respectively, examples of stable, periodic and chaotic behaviours.

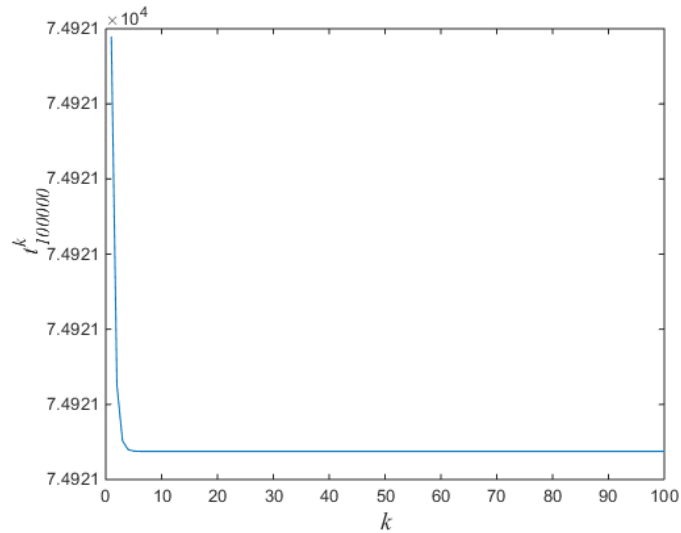

Figure 5 – Stable solution of (11) for $n = 100{,}000$ and $\delta = 0.1595388$.

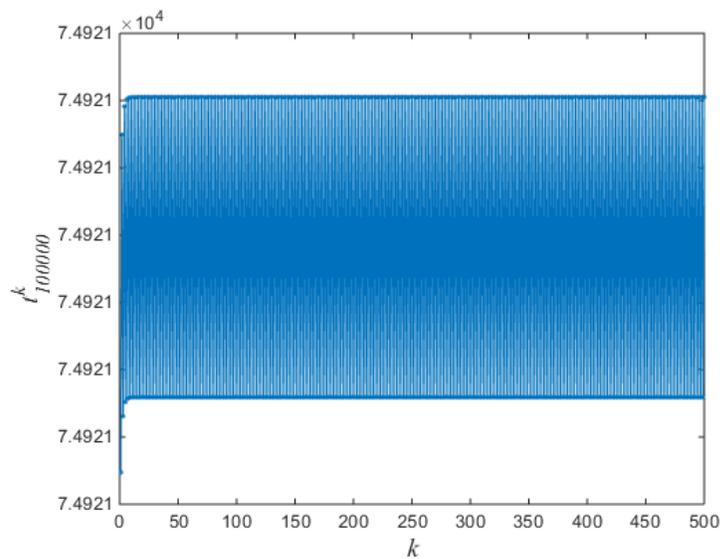

Figure 6 – Periodic solution of (11) for $n = 100{,}000$ and $\delta = 0.5$.

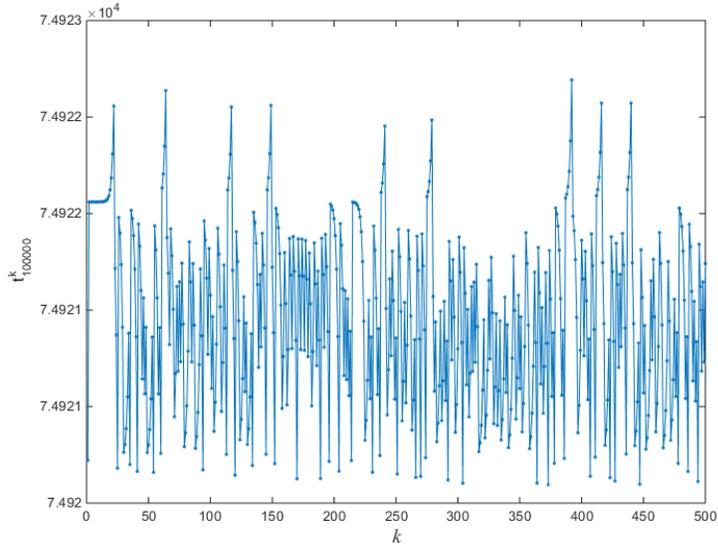

Figure 7 – Chaotic solution of (11) for $n = 100{,}000$ and $\delta = 2$.

On the other hand, Fig. 8 shows the bifurcation diagram when $n = 100{,}000$.

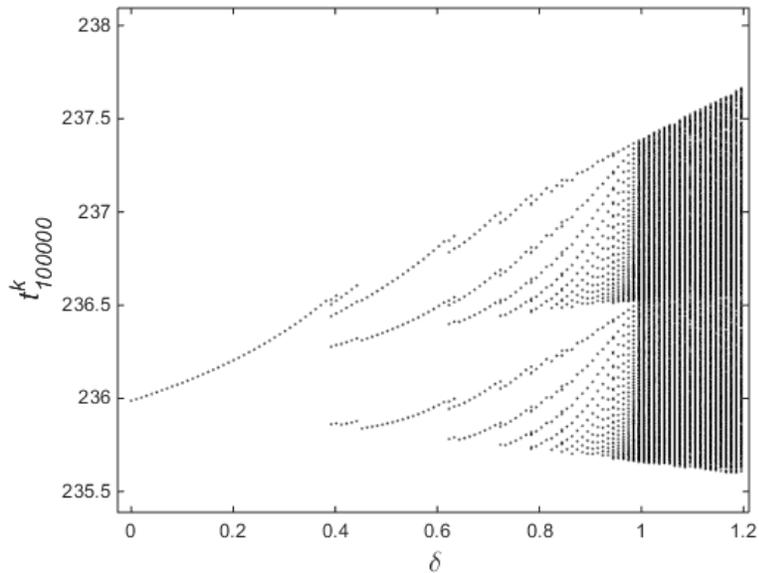

Figure 8 – Bifurcation diagram of (11) versus $\delta$ for $n = 100{,}000$.

## 5. Conclusions

It was provided a non-linear difference equation that can converge to the zeros of the Riemann zeta-function on the critical line. The zeros obtained are good enough to test properties of the zeta-function or to calculate functions that depend on the zeros, like the discrete cosine Riemann transform that can be used for integer factorization. The non-linear map proposed has a rich dynamic with stable, periodic and chaotic solutions, depending on the value of $\delta$. Obviously, the zeros on the critical line represent stable solutions and the value of $\delta$ that maximizes the accuracy of $t_n$ calculated using (11) is the value for which the first bifurcation occurs.

**Acknowledgements**

This study was financed in part by the Coordenação de Aperfeiçoamento de Pessoal de Nível Superior - Brasil (CAPES) - Finance Code 001, and CNPq via Grant no. 307062/2014-7. Also, this work was performed as part of the Brazilian National Institute of Science and Technology for Quantum Information.